\documentclass{conm-p-l} 
\usepackage{xypic}
\xyoption{ps}\xyoption{dvips} 
\newbox\mybox 
\def\overtag#1#2#3{\setbox\mybox\hbox{$#1$}\hbox to
  0pt{\vbox to 0pt{\vglue-#3\vglue-\ht\mybox\hbox to \wd\mybox
      {\hss$\ss#2$\hss}\vss}\hss}\box\mybox}
\def\undertag#1#2#3{\setbox\mybox\hbox{$#1$}\hbox to 0pt{\vbox to
    0pt{\vglue#3\vglue\ht\mybox\hbox to \wd\mybox
      {\hss$\ss#2$\hss}\vss}\hss}\box\mybox} 
\def\lefttag#1#2#3{\hbox to 0pt{\vbox to 0pt{\vss\hbox to
      0pt{\hss$\ss#2$\hskip#3}\vss}}#1} 
\def\righttag#1#2#3{\hbox to 0pt{\vbox to 0pt{\vss\hbox to
      0pt{\hskip#3$\ss#2$\hss}\vss}}#1} 
\let\ss\scriptstyle
\def\notags{\def\overtag##1##2##3{##1}
  \def\undertag##1##2##3{##1}\def\lefttag##1##2##3{##1}
  \def\righttag##1##2##3{##1}}
\def\Dot{\lower.2pc\hbox to 2pt{\hss$\bullet$\hss}}
\def\Circ{\lower.2pc\hbox to 2pt{\hss$\circ$\hss}}
\def\Vdots{\raise5pt\hbox{$\vdots$}}
\newcommand\lineto{\ar@{-}}
\newcommand\dashto{\ar@{--}}
\newcommand\dotto{\ar@{.}}
\def\WDN#1.{\marginpar{\rightskip=0pt plus2\hsize#1}} 
\def\WDN#1.{}                         

\newcommand\Aut{\operatorname{Aut}}
\newcommand\splice{\Gamma}
\newcommand\rsplice{\Gamma_\bullet}

\renewcommand\vert{\operatorname{vert}}

\renewcommand\deg{\operatorname{deg}}

\newcommand\R{{\mathbb R}}
\newcommand\C{{\mathbb C}}

\newtheorem{theorem}{Theorem}[section]

\newtheorem{proposition}[theorem]{Proposition}

\theoremstyle{definition}
\newtheorem{definition}[theorem]{Definition}
\newtheorem{problem}[theorem]{Problem}

\begin{document}
\title{Algorithms for polynomials in two variables}
\author{Walter D. Neumann}
\address{Department of Mathematics and Statistics\\The University of
Melbourne\\Parkville, Vic 3052\\Australia}
\email{neumann@ms.unimelb.edu.au}
\author{Penelope G. Wightwick}
\address{Department of Mathematics and Statistics\\The University of
Melbourne\\Parkville, Vic 3052\\Australia}
\email{pw@ms.unimelb.edu.au}
\subjclass{12Y05, 14H50}
\keywords{polynomial automorphism}
\thanks{This research was supported by the Australian Research Council
and an Australian Postgraduate Research Award.}
\begin{abstract}
  Vladimir Shpilrain and Jie-Tai Yu have asked for an effective
  algorithm to decide if two elements of $\C[x,y]$ are related by an
  automorphism of $\C[x,y]$.  We describe here an efficient algorithm,
  due to the second author, that decides this question and finds the
  automorphism if it exists.  We also discuss some examples related to
  work of Shpilrain and Yu. Part of the purpose of this paper is to
  advertise the use of splice diagrams in studying $\C[x,y]$.
\end{abstract} 
\maketitle

\section{Introduction}
\label{sec:intro}

We give an efficient algorithmic answer to the following problem,
which was posed by Shpilrain and Yu \cite{shpilrain-yu-99} and resolved
by them in special cases. The algorithm is due to the second author in
\cite{wightwick}, where it is described in terms of Newton polygons.
Here we describe it in terms of splice diagrams.
\begin{problem}\label{qu:1}
  Given polynomials $f,g\in\C[x,y]$, decide if
  there exists a polynomial automorphism $\phi\in\Aut\C[x,y]$ with
  $f\phi=g$ and find $\phi$ if it exists.
\end{problem}
As we discuss in the final section, $\C$ could be replaced by any
field, although in finite characteristic the description in terms of
splice diagrams given here would need some modification.

Any $\phi\in\Aut\C[x,y]$ induces a polynomial bijection $\C^2\to\C^2$
(which we also call $\phi$) and vice versa.  From this point of view,
the problem asks for a polynomial bijection $\phi\colon\C^2\to\C^2$
which makes the diagram 
\begin{equation*}
  \xymatrix{
\C^2\ar[r]^{\phi}\ar[d]^g&\C^2\ar[d]^f \\
\C\ar@{=}[r]&\C
}
\end{equation*}
commutative. We shall use the terminology ``$f$ is right-equivalent to
$g$'' for this.  From a topological point of view it is also natural
to ask about ``right-left equivalence'', which asks for polynomial
bijections $\phi\colon\C^2\to\C^2$ and $\psi\colon\C\to\C$ making
\begin{equation*}
  \xymatrix{
\C^2\ar[r]^{\phi}\ar[d]^g&\C^2\ar[d]^f \\
\C\ar[r]^\psi&\C
}
\end{equation*}
commutative ($\psi$ is an affine map $\psi(z)=az+b$ with $a\ne0$). Our
results apply also to this. We shall use the abbreviation
``equivalence'' to mean right-equivalence.

Shpilrain and Yu say $f$ and $g$ in $\C[x,y]$ are ``isomorphic'' if
the curves $f(x,y)=0$ and $g(x,y)=0$ are isomorphic as affine
schemes; that is, their rings of functions $\C[x,y]/(f)$ and
$\C[x,y]/(g)$ are isomorphic.  In the final section we discuss some
of their examples, and additional examples motivated from their work,
of isomorphic but non-equivalent polynomials.

We first recall the structure of the automorphism group $\Aut\C[x,y]$.
It has two subgroups
\begin{align*}
  A&:=\operatorname{Aff}\C^2=\{(x,y)\mapsto(ax+by+s,cx+dy+t):
a,b,c,d,s,t\in\C\}\\
  B&:=\{(x,y)\mapsto(ax+g(y),dy+t):g(y)\in\C[y], a,d,t\in\C\}
\end{align*}
$B$ is called  
the ``Jonqui\`ere subgroup''.
\begin{theorem}[Jung \cite{jung}, 1942] \label{th:jung}
  $\Aut\C[x,y]$ is generated by its subgroups $A$ and $B$.
\end{theorem}
This theorem implies that any $\phi\in\Aut\C[x,y]$ has an expression
\begin{equation*}
  \phi=\phi_1\phi_2\dots\phi_n
\end{equation*}
with $\phi_i$ in $A$ for $i$ even and in $B$ for $i$ odd or vice
versa. Moreover, we may assume that if $n>1$ then no $\phi_i$ is in
$A\cap B$.
We call such an expression \emph{normal form}.

If each $\phi_i$ in $B$ is actually in the subgroup of
\emph{triangular automorphisms}
$$
T:=\{(x,y)\mapsto(x+g(y),y):g(y)\in\C[y]\}.
$$
we shall speak of \emph{strict normal form}. We can always change a
normal form representation with $n>1$ to strict normal form by
applying changes of the form: replace $\phi_i$ and $\phi_{i+1}$ by
$\phi_i\chi$ and $\chi^{-1}\phi_{i+1}$ with $\chi\in A\cap B$.
  
The following theorem is stated in slightly different form in
\cite{wightwick}.
\begin{theorem}[Wightwick \cite{wightwick}]\label{th:pw}
  Suppose $f\in\C[x,y]$ is a non-constant polynomial and
  $\phi\in\Aut\C[x,y]$. Suppose $\phi=\phi_1\phi_2\dots\phi_n$ is a
  normal form expression for $\phi$, and put
  $f_i=f\phi_1\phi_2\dots\phi_i$ for $i=1,\dots,n$, and $f_0=f$.  Then
  the sequence of degrees $\deg(f_i)$ satisfies
\begin{equation*}
  \deg(f_0)\ge\dots\ge\deg(f_k)\le\dots\le\deg(f_n),\quad\text{for
  some }0\le k\le n.
\end{equation*}
Moreover, there is at most one $\phi_i\in B$ for which
$\deg(f_{i-1})=\deg(f_i)$, and if this happens then $\deg(f_i)$ is the
minimal degree $\deg(f_k)$.
\end{theorem}

This theorem implies the result (van der Kulk \cite{kulk}, 1953) that
$\Aut\C[x,y]$ is the amalgamated free product of $A$ and $B$,
amalgamated along $A\cap B$. For if not, there would be a non-trivial
normal form representation of $\phi = 1$ with $\phi_1\in B$, and
applying the theorem with this $\phi$ to $f=x$ gives a contradiction.

The theorem clearly implies that $n$ is at most
$2(\deg(f_0)+\deg(f_n))+1$.  It is also easy to see that that each
$\phi_i$ has degree at most $\max(\deg(f_{i-1}),\deg(f_i))$ (this also
follows from the proof in \cite{wightwick}).  This already implies an
algorithm for Problem \ref{qu:1}, since it implies a bound of
$N^{4N+1}$ on the degree of the automorphism $\phi$ with
$N=\max(\deg(f), \deg(g))$, so finding $\phi$ amounts to solving a
system of algebraic equations in the coefficients of this
automorphism, for which computational techniques are known.  However,
we promised an efficient algorithm, and this algorithm would be
inefficient even if this bound on $\deg(\phi)$ were much improved (as
it can be).

We will describe the efficient algorithm in terms of the splice
diagrams introduced in the book \cite{eisenbud-neumann}. They were
used there to describe the local topology of plane curve
singularities. Later, in \cite{neumann1}, they were used to discuss the
global topology of polynomial maps $f\colon\C^2\to\C$, which is what
interests us here. 

\section{Splice diagrams}

Associated to a fiber $f^{-1}(c)\subset\C^2$ of a non-constant
polynomial $f\in\C[x,y]$ is a combinatorial invariant: the
\emph{
splice diagram} $\splice(f^{-1}(c))$ for $f^{-1}(c)$.
It is a tree with certain numerical and other decorations on it.
There are finitely many values of $c$ for which this diagram differs
from its generic value.  The generic diagram is called the
\emph{regular splice diagram} for $f$, denoted $\splice(f)$. As we
shall describe, $\splice(f^{-1}(c))$ encodes behaviour at infinity of
$f^{-1}(c)$.  Therefore, for the finitely many $c$ for which
$\splice(f^{-1}(c))$ is not regular we speak of $f^{-1}(c)$ being
\emph{irregular at infinity} or having ``singularities at infinity.''

For example, the so called Brian\c con polynomial,
$$f(x,y)=x^2(1+xy)^4+3x(1+xy)^3+(3-\frac83x)(1+xy)^2 -4(1+xy)+y,$$
discussed in detail in \cite{artal-cassou-luengo}, has regular splice
diagram: 

$$
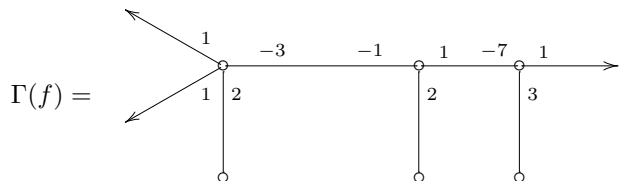
\notags
\objectmargin{0pt}\spreaddiagramrows{-5pt}
\spreaddiagramcolumns{12pt}
\splice(f)=\quad\diagram
\\
&\overtag\Circ{(0)}{16pt}
\ulto_(.25){1}\dlto^(.25){1}
\ddline^(.25){2}\rrline^(.25){-3}
^(.75){-1}&&\overtag\Circ{(2)}{16pt}
\ddline^(.25){2}\rline^(.25){1}^(.75){-7}
&\overtag\Circ{(3)}{16pt}\ddline^(.25){3}\rto^(.25){1}&\\ \\
&\lefttag\Circ{(0)}{6pt}&&
\lefttag\Circ{(1)}{6pt}&\lefttag\Circ{(1)}{6pt}
\enddiagram
$$
and two fibers which are irregular at infinity:
$$\notags
\objectmargin{0pt}\spreaddiagramrows{-5pt}
\spreaddiagramcolumns{12pt} \splice(f^{-1}(0))=\quad
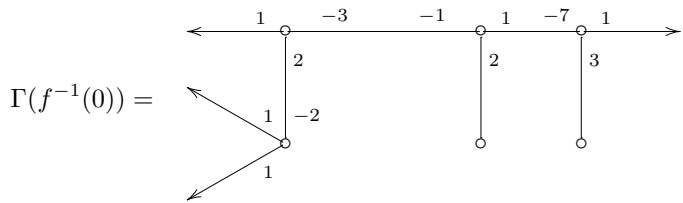
\diagram &\overtag\Circ{(0)}{16pt}
\lto_(.25){1}\ddline^(.25){2}^(.75){-2}
\rrline^(.25){-3}
^(.75){-1}&&\overtag\Circ{(2)}{16pt}
\ddline^(.25){2}\rline^(.25){1}^(.75){-7}
&\overtag\Circ{(3)}{16pt}\ddline^(.25){3}\rto^(.25){1}&\\ \\
&\righttag\Circ{(-1)}{6pt}\ulto_(.25){1}\dlto^(.25){1}&&
\lefttag\Circ{(1)}{6pt}&
\lefttag\Circ{(1)}{6pt}\\
& \enddiagram
$$
$$\notags
\objectmargin{0pt}\spreaddiagramrows{-5pt}
\spreaddiagramcolumns{12pt}
\splice(f^{-1}({\scriptstyle\frac{-16}{9}}))=\quad
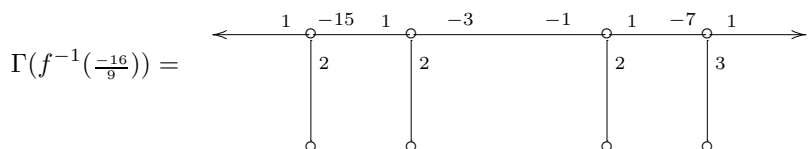
\diagram\\
&\lto_(.25){1}\overtag\Circ{(-6)}{16pt}
\rline^(.25){-15}^(.75){1}\ddline^(.25){2}
&\overtag\Circ{(0)}{16pt}\ddline^(.25){2}\rrline^(.25){-3}
^(.75){-1}&&\overtag\Circ{(2)}{16pt}
\ddline^(.25){2}\rline^(.25){1}^(.75){-7}
&\overtag\Circ{(3)}{16pt}\ddline^(.25){3}\rto^(.25){1}&\\ \\
&\lefttag\Circ{(-3)}{6pt}&\lefttag\Circ{(0)}{6pt}&&
\lefttag\Circ{(1)}{6pt}
&\lefttag\Circ{(1)}{6pt}
\enddiagram
$$

We explain such diagrams in greater detail below, but first we explain
why they are useful invariants.

First, $\splice(f^{-1}(c))$ is easily computable (by hand; or explicit
computer code is also available in Magma \cite{magma}).  In fact,
$\splice(f^{-1}(c))$ is simply a graphical representation of the
exponents of the topologically relevant terms in the Puiseux
expansions at infinity of $f(x,y)=c$, and Newton already knew how to
compute Puiseux expansions. The regular splice diagram $\splice(f)$ is
easily derived purely combinatorially from $\splice(f^{-1}(c))$ (see
\cite{neumann2}), so one does not need to know a regular value of $f$
to find $\splice(f)$.

Second, each $\splice(f^{-1}(c))$ individually is an invariant under
equivalence of $f$, while the collection of them is an invariant under
right-left equivalence (in particular, the regular splice diagram
$\splice(f)$ is invariant under right-left equivalence of $f$).

Actually, splice diagrams exist in two versions, the rooted and
unrooted splice diagrams, and the above are the unrooted versions.
The rooted splice diagrams will be denoted
$\rsplice(f^{-1}(c))$ and $\rsplice(f)$ respectively. For the above
examples the rooted diagrams differ from the unrooted ones by the
addition of a root vertex (drawn as a filled dot ``$\bullet$'' with
adjacent numeric weights $1$) in the center of the edge with weights
$-3,-1$.

We usually omit edge weights $1$. For example, the rooted
splice diagram for Brian\c con's polynomial is therefore
$$
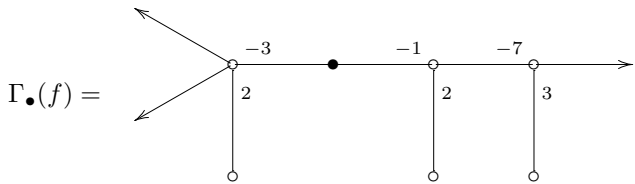
\notags
\objectmargin{0pt}\spreaddiagramrows{-5pt}
\spreaddiagramcolumns{12pt}
\rsplice(f)=\quad\diagram
\\
&\overtag\Circ{(0)}{16pt}
\ulto\dlto\ddline^(.25){2}\rline^(.25){-3}
&\Dot\rline
^(.75){-1}&\overtag\Circ{(2)}{16pt}
\ddline^(.25){2}\rline
^(.75){-7}
&\overtag\Circ{(3)}{16pt}\ddline^(.25){3}\rto&\\ \\
&\lefttag\Circ{(0)}{6pt}&&
\lefttag\Circ{(1)}{6pt}&\lefttag\Circ{(1)}{6pt}
\enddiagram
$$

The unrooted splice diagrams encode topology of $f\colon\C^2\to\C$,
while the rooted splice diagrams encode the same topology together
with its relationship to a generic line in $\C^2$.  Thus, the rooted
splice diagram is \emph{not} an invariant of equivalence, since the
concept of generic line is not invariant under polynomial
automorphisms.  The topological meaning of these diagrams is detailed
in \cite{neumann1} and surveyed in \cite{neumann-proc}.  We will here
give a brief description of the algebraic meaning. See also
\cite{cassou} for an extended survey and \cite{wightwick-msc} for
details of how to compute splice diagrams using Newton polygons.

We first describe the general structure of the splice diagrams that
arise in our situation (the splice diagrams used to describe plane
curve singularities are similar except that edge determinants are
positive in Item \ref{it:rpi2}).

\subsection{Properties of splice diagrams}
\subsubsection{} Each \emph{arrowhead} in the splice diagram 
corresponds to a branch at infinity of $f^{-1}(c)$ (a \emph{branch at
  infinity} can be taken to be a component of the intersection of
$f^{-1}(c)$ with the complement of a large ball in $\C^2$).
\subsubsection{} Each edge emanating from the \emph{root vertex} 
corresponds to a point at infinity of $f^{-1}(c)$ (a \emph{point at
  infinity} is a point of intersection of the closure of $f^{-1}(c)$
in $\C P^2$ with the line at infinity $\C P^1$). All edge weights
adjacent to the root vertex are $1$.
\subsubsection{} Any other vertex of the splice diagram is either
\begin{itemize}
\item a \emph{leaf}: a non-arrowhead vertex of valency 1. It has no
  adjacent weights.
\item a \emph{node}: a non-arrowhead vertex of valency $\ge2$. It has
  pairwise coprime integer weights associated to its adjacent edges.
\end{itemize}
Each edge will thus have a weight at its near end, seen from the node,
and may have a weight at its far end.  We call weights ``near'' or
``far'' correspondingly.
\subsubsection{}\label{it:rpi1} All near weights are positive and at 
most one of the near weights at any node differs from $1$ (edge
weights 1 are usually omitted in the diagram). A far weight can be any
integer.
\subsubsection{}\label{it:rpi2} The \emph{edge determinant} of an 
edge connecting two nodes (one of which may be the root vertex) is the
product of the two weights on the edge minus the product of the
weights adjacent to the edge. All edge determinants are negative. For
example, in the rooted diagram above the three edge determinants are,
from left to right, $-3-2=-5$, $-1-2=-3$, and
$-7-(-1\times2\times3)=-1$.
\subsubsection{} ({\bf Reduction 1}) Any edge with edge weight $1$ 
leading to a leaf may be deleted (with the leaf).\label{red1}
\subsubsection{} ({\bf Reduction 2}) Any node (other than the 
root vertex) of valency $2$ may be removed, and the adjacent edges
coalesced into a single edge.\label{red2}

\begin{definition}
\label{de:reduced}
A diagram is \emph{reduced} if neither Reduction 1 nor Reduction 2 can
be applied to it. In particular, all non-root nodes have valency
$\ge3$.

We will assume diagrams are reduced unless otherwise stated, and the
notations $\rsplice(f^{-1}(c))$, $\rsplice(f^{-1}(c))$,
$\splice(f^{-1}(c))$, and $\splice(f)$ will always mean reduced
diagrams.

Here, the \emph{unrooted splice diagram} $\splice(f^{-1}(c))$ or
$\splice(f)$ is the result of changing the root vertex of the
corresponding rooted diagram to an ordinary vertex and then reducing
by Reductions 1 and/or 2 as necessary. However, there is an
exceptional case described at the end of this section, where a further
reduction must be done to eliminate a vanishing edge determinant.
\end{definition}

\subsection{Splice diagrams via Puiseux expansion}
The rooted diagram relates to the algebra as follows. Pick an
arrowhead and consider the corresponding branch at infinity of
$f(x)=c$. Change coordinates linearly so this branch occurs at the
point at infinity $[x{:}y]=[1{:}0]$ and consider the Puiseux expansion
at this branch in the ``Newton form''
$$y=x^{q_1/p_1}(a_1+x^{d_1/p_1p_2}(a_2+x^{d_2/p_1p_2p_3}+\dots)\dots)).$$
Here $q_1<p_1$ and $d_i<0$ for each $i$.  One can show that there is a
$k$ such that $p_i=1$ for $i> k$.  We only need the terms up to this
point.  We then form the (maybe unreduced) splice diagram
$$\notags
\objectmargin{0pt}\spreaddiagramrows{-5pt}
\spreaddiagramcolumns{12pt}
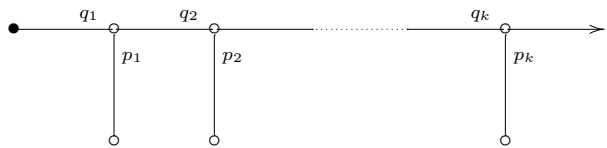
\diagram
\Dot\rline^(.75){q_1}&\Circ\ddline^(.25){p_1}\rline^(.75){q_2}&\Circ
\ddline^(.25){p_2}\rline&\dotto[r]
&\rline^(.75){q_k}&\Circ\ddline^(.25){p_k}\rto&\\ \\
&\Circ&\Circ&&&\Circ
\enddiagram
$$
where the $q_i$ for $i\ge2$ are chosen inductively so that $d_i$ is
the edge determinant $q_{i+1}-q_ip_{i+1}$ for each $i$.  Finally, we
assemble these diagrams for the different branches at infinity of
$f(x)=c$ by merging the parts that correspond to identical initial
segments of the Puiseux expansions and then reduce the diagram as
necessary.

For example, the right-hand arrowhead in the above diagrams for the
Brian\c con polynomial corresponds to a branch of $f(x,y)=c$ where the
Puiseux expansion at infinity is (writing the expansion as $x$ in
terms of $y$ since the branch is at $[0{:}1]$ rather than $[1{:}0]$)
$$x=y^{-1/2}(i+y^{-1/6}(\frac i2+\dots)\dots).$$
This gives the diagram
$$
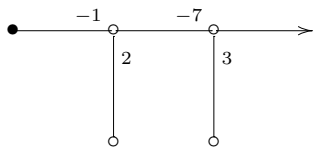
\notags
\objectmargin{0pt}\spreaddiagramrows{-5pt}
\spreaddiagramcolumns{12pt}
\diagram
\Dot\rline
^(.75){-1}&\overtag\Circ{(2)}{16pt}
\ddline^(.25){2}\rline
^(.75){-7}
&\overtag\Circ{(3)}{16pt}\ddline^(.25){3}\rto&\\ \\
&\lefttag\Circ{(1)}{6pt}&\lefttag\Circ{(1)}{6pt}
\enddiagram
$$

There are generically two branches of $f(x,y)=c$ at the point
$[1{:}0]$ at infinity and they have Puiseux expansions
$$y=x^{-1}(-1+x^{-1/2}(\alpha+\dots)\dots)$$
with
$\alpha=\sqrt{4/3\pm\sqrt{16/9+c}}$. For $c\ne0,-16/9$ this gives two
branches each with (unreduced) diagram
$$\notags \objectmargin{0pt}\spreaddiagramrows{-5pt}
\spreaddiagramcolumns{12pt} 
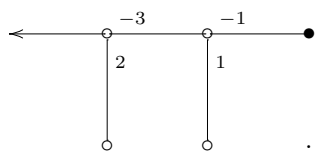
\diagram
&\lto\Circ\ddline^(.25){2}\rline^(.25){-3}&
\Circ\rline^(.25){-1}\ddline
^(.25){1}&\Dot\\ \\
&\Circ&\Circ&. \enddiagram
$$
Combining these two branches gives the (unreduced)
diagram
$$\notags \objectmargin{0pt}\spreaddiagramrows{-5pt}
\spreaddiagramcolumns{12pt}
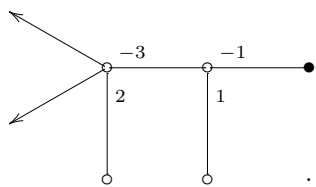
\diagram\\
&\ulto\dlto\Circ\ddline^(.25){2}\rline^(.25){-3}&
\Circ\rline^(.25){-1}\ddline
^(.25){1}&\Dot\\ \\
&\Circ&\Circ&.  \enddiagram
$$
which reduces by reduction moves of Items \ref{red1}.\ and
\ref{red2}.\ above to 
$$\notags
\objectmargin{0pt}\spreaddiagramrows{-5pt}
\spreaddiagramcolumns{12pt}
\diagram\\
&\ulto\dlto\Circ\ddline^(.25){2}\rline^(.25){-3}&\Dot\\ \\
&\Circ&.
\enddiagram
$$
Combining this with the diagram for the branch at $[0{:}1]$ gives the
regular rooted splice diagram shown earlier.

The two values $c=0,-16/9$ are clearly special for the above Puiseux
expansion. When $c=0$, in addition to the branch at $[1{:}0]$
described above with $\alpha=\sqrt{8/3}$, there are two branches
corresponding to $\alpha=0$:
$$
y=x^{-1}(-1+x^{-1}(\beta+\dots)\dots),\quad\beta=-\frac34\pm\sqrt3.$$
When $c=-16/9$, so $\alpha=\sqrt{4/3}$, there is just one branch at
$[1{:}0]$ and it has an additional relevant term in its Puiseux
expansion:
$$y=x^{-1}(-1+x^{-1/2}(\sqrt{4/3}+x^{-1/2}(-{
3/4}+\dots)\dots)).$$
These give the two irregular diagrams shown earlier.

\subsection{Exceptional splice diagrams}\label{subsec:exceptional}
As mentioned earlier in Definition \ref{de:reduced}, there is an
exceptional case for the relationship between reduced rooted and
unrooted diagrams. The following rooted splice diagram will be called
\emph{exceptional} if $p\ge q\ge1$, $gcd(p,q)=1$, and $pd>1$:
$$\def\dec{\ar[d]+<10pt,0pt>\ar[d]-<10pt,0pt>\ar@{}[d]*{\dots}}
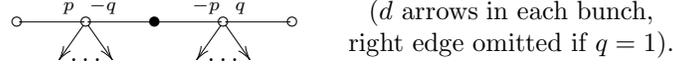
\xymatrix@R=12pt@C=24pt@M=0pt@W=0pt@H=0pt{
  \Circ\lineto[r]^(.75){p}&\Circ\lineto[r]^(.25){-q}\dec&
  \Dot\lineto[r]^(.75){-p}&\Circ\lineto[r]^(.25){q}\dec&\Circ\\
  &&&&\\
  }
\quad \txt<12pc>{($d$ arrows in each bunch,\\
  right edge omitted if $q=1$).}$$
This is the regular reduced rooted splice diagram for
$f(x,y)=x^{dp}y^{dq}$, for example.

(This $f(x,y)$ has one irregular fiber $f^{-1}(0)$ with splice diagram
$$\xymatrix@M=0pt@W=0pt{(dp)~&{\bullet}\ar[l]\ar[r]&~(dq)}.$$
The
weights in parentheses refer to multiplicity of nonreduced branches at
infinity. There are other polynomials with the same regular rooted
splice diagram which have different irregular ones, see
\cite{neumann1}.)

The reduced unrooted splice diagram for this rooted diagram is
$$
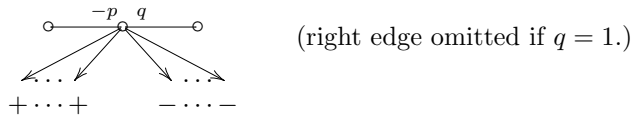
\def\decl{\ar[dddl]+<10pt,0pt>\ar[dddl]-<10pt,0pt>\ar@{}[dddl]*{\dots}}
\def\decr{\ar[dddr]+<10pt,0pt>\ar[dddr]-<10pt,0pt>\ar@{}[dddr]*{\dots}}
\xymatrix@R=6pt@C=12pt@M=0pt@W=0pt@H=0pt{\Circ\lineto[r]^(.75){-p}
&\Circ\lineto[r]^(.25){q}
\decl
\decr
&\Circ\\ \\ \\&&&\\
+\cdots+&&-\cdots-
}\quad\txt{~\\(right edge omitted if $q=1$.)}
$$
There are $d$ arrows with $+$ signs and $d$ with $-$ signs. The
signs refer to orientations. See \cite{neumann1} for an explanation.

\subsection{Non-reduced fibers} 
As above, if a polynomial $f$ has non-reduced fibers then weights at
arrowheads of the splice diagram are used to indicate multiplicities
of components of such a fiber.  Such a diagram is of course an
irregular diagram.

\section{Invariants from splice diagrams}\label{sec:invariants}

Much of this section is not essential for what follows, but the
invariants encoded in splice diagrams are often helpful in computing
them, since they give strong numerical constraints.

Given a vertex $v$ of a splice diagram $\splice$, the \emph{linking
  coefficient at $v$} (also called \emph{multiplicity at $v$}) is the
sum over arrowheads $w$ of $\splice$:
$$\ell_v=\sum_w\ell(v,w),$$
where $\ell(v,w)$ is the product of edge
weights directly adjacent to but not on the path from $v$ to $w$ in
$\splice$.  For example, in the regular and irregular
rooted diagrams for the Brian\c
con polynomial these linking coefficients are as indicated in
parentheses in the diagrams:
$$
\objectmargin{0pt}\spreaddiagramrows{-5pt}
\spreaddiagramcolumns{12pt}
\rsplice(f)=\quad\diagram
\\
&\overtag\Circ{(0)}{14.1pt}
\ulto\dlto\ddline^(.25){2}\rline^(.25){-3}
&\overtag\Dot{(10)}{14.1pt}\rline
^(.75){-1}&\overtag\Circ{(2)}{14.1pt}
\ddline^(.25){2}\rline
^(.75){-7}
&\overtag\Circ{(3)}{14.1pt}\ddline^(.25){3}\rto&\\ \\
&\lefttag\Circ{(0)}{6pt}&&
\lefttag\Circ{(1)}{6pt}&\lefttag\Circ{(1)}{6pt}
\enddiagram
$$
$$
\objectmargin{0pt}\spreaddiagramrows{-5pt}
\spreaddiagramcolumns{12pt} \rsplice(f^{-1}(0))=\quad
\diagram &\overtag\Circ{(0)}{14.1pt}
\lto_(.25){1}\ddline^(.25){2}^(.75){-2}
\rrline^(.25){-3}&\overtag\Dot{(10)}{14.1pt}\rline
^(.75){-1}&\overtag\Circ{(2)}{14.1pt}
\ddline^(.25){2}\rline^(.25){1}^(.75){-7}
&\overtag\Circ{(3)}{14.1pt}\ddline^(.25){3}\rto^(.25){1}&\\ \\
&\righttag\Circ{(-1)}{6pt}\ulto_(.25){1}\dlto^(.25){1}&&
\lefttag\Circ{(1)}{6pt}&
\lefttag\Circ{(1)}{6pt}\\
& \enddiagram
$$
$$
\objectmargin{0pt}\spreaddiagramrows{-5pt}
\spreaddiagramcolumns{12pt}
\rsplice(f^{-1}({\scriptstyle\frac{-16}{9}}))=\quad\diagram\\
&\lto_(.25){1}\overtag\Circ{(-6)}{14.1pt}
\rline^(.25){-15}^(.75){1}\ddline^(.25){2}
&\overtag\Circ{(0)}{14.1pt}\ddline^(.25){2}\rrline^(.25){-3}
&\overtag\Dot{(10)}{14.1pt}\rline
^(.75){-1}&\overtag\Circ{(2)}{14.1pt}
\ddline^(.25){2}\rline^(.25){1}^(.75){-7}
&\overtag\Circ{(3)}{14.1pt}\ddline^(.25){3}\rto^(.25){1}&\\ \\
&\lefttag\Circ{(-3)}{6pt}&\lefttag\Circ{(0)}{6pt}&&
\lefttag\Circ{(1)}{6pt}
&\lefttag\Circ{(1)}{6pt}
\enddiagram
$$

\subsection{Properties relating to linking coefficients}
\subsubsection{} \label{it:degree} The degree of a polynomial is 
always the linking coefficient $\ell_{\bullet}$ at the root vertex.
\subsubsection{} A splice diagram is regular if and only if it 
has no negative linking coefficients and the fiber in question is
reduced (i.e., there are no multiplicity weights at arrowheads).
\subsubsection{} A polynomial has fibers with irregular splice 
diagrams if and only if the regular splice diagram has at least one
zero linking coefficient.  Each irregular diagram determines the
regular one and the regular one strongly constrains the number and
form of the irregular ones (\cite{neumann2}).
\subsubsection{} \label{it:reg}The Euler characteristic of the generic 
fiber of $f$ is
$$\chi_{reg}:=\sum_{v\in\vert\splice(f)}(2-\delta_v)\ell_v,$$
where
$\vert\splice(f)$ is the set of non-arrowhead vertices of $\splice(f)$
and $\delta_v$ is valency of vertex $v$ (number of edges at $v$).
\subsubsection{} \label{it:irr}For any $c\in\C$ define the 
\emph{Milnor number at infinity} of $f^{-1}(c)$ as
$$\lambda_c:=\sum_{\genfrac{}{}{0pt}{}{\scriptstyle
    v\in\vert\splice(f^{-1}(c))}{\scriptstyle\ell_v<0}}
(2-\delta_v)\ell_v,$$
so $\lambda_c=0$ unless $\splice(f^{-1}(c))$ is
irregular.  Moreover, if $f^{-1}(c)$ is a reduced fiber,
define its \emph{total Milnor number} as the sum of Milnor numbers:
$$\mu_c:=\sum_{p\in f^{-1}(p)}\mu_p$$
(since $\mu_p=0$ at a
non-singular point, this is a finite sum). Then the Euler
characteristic of a reduced fiber $f^{-1}(c)$ is given by
$$\chi(f^{-1}(c))=\mu_c+\sum_{v\in\vert\splice(f^{-1}(c))}
(2-\delta_v)\ell_v=\mu_c+\lambda_c+\chi_{reg}.$$
\subsubsection{} (Suzuki \cite{suzuki})\label{it:suz} For any $f$ one 
has $1-\chi_{reg}=\sum_{c\in\C}\bigl(\chi(f^{-1}(c)-\chi_{reg}\bigr)$,
so if all fibers of $f$ are reduced (i.e., $f$ has only isolated
singularities) then
  $$1-\chi_{reg}=\sum_{c\in\C}(\mu_c+\lambda_c).$$

\bigbreak
In the formulae of Items \ref{it:reg}. and \ref{it:irr}. one can use
either the rooted or unrooted diagram (in fact the the formulae do not
need reduced diagrams).  

For any polynomial, as soon as singularities and ``singularities at
infinity'' are found contributing sufficient $\mu_c$ and $\lambda_c$
to satisfy Suzuki's formula, one knows that all non-generic fibers
have been found.  For example, for Brian\c con's polynomial, Items
\ref{it:reg}. and \ref{it:irr}. give:
$$\chi_{reg}=-3;\quad\lambda_0=1,\quad \lambda_{-16/3}=3.$$
These
already satisfy Suzuki's formula, so the fibers $f^{-1}(0)$ and
$f^{-1}(-16/3)$ are indeed the only fibers with singularities at
infinity and there are no ``finite'' singularities.  Brian\c con's
polynomial was the first non-singular polynomial discovered that also
has all fibers connected \cite{artal-cassou-luengo}, although Brian\c
con constructed it for a different purpose.

\section{Newton polygons and splice diagrams}\label{sec:newton}

Given a non-constant polynomial $f(x,y)$, our first step will be to
put it into as simple form as possible. We first describe how the part
of $\rsplice(f)$ nearest the root can be read off from the Newton
polygon of $f$. For more details see \cite{neumann1} or
\cite{wightwick-msc}.

The \emph{Newton polygon} of $f$ is the convex hull of the set of
points $$\{(a,b)\in \R^2:x^ay^b\text{ is a monomial of
}f\}\cup\{(0,0)\}.$$ (This would better be called the ``regular Newton
polygon'' and the adjective ``regular'' dropped if $(0,0)$ is not
explicitly added, but we shall only use this version.)

If the Newton polygon consists of a single line segment then either:
\begin{itemize}
\item $f$ is a polynomial in just one of its variables,
in which case the regular rooted splice diagram is 
$$
\left.\vcenter{
    \xymatrix@R=4pt@C=18pt@M=0pt@W=0pt@W=0pt{& & \\
      \Dot\lineto[r]^(.75){0}&\Circ\ar[ur]\ar[dr]&\Vdots\\
      & &}}~~\right\}\operatorname{deg}(f) \quad\text{if }
\operatorname{deg}(f)>1,\qquad
\xymatrix@R=4pt@C=18pt@M=0pt@W=0pt@W=0pt{ \Dot\ar[r]&} \quad
\text{if }\operatorname{deg}(f)=1,$$
\item $f(x,y)=g(x^py^q)$ for some polynomial
$g$ of degree $d$, say, in which case the regular rooted splice
diagram is 
$\xymatrix@1@R=4pt@C=18pt@M=0pt@W=0pt@W=0pt{&\ar[l]\bullet\ar[r]&}$
if $p=q=d=1$, and of the exceptional type of Subsection
\ref{subsec:exceptional} otherwise.
\end{itemize}

We assume from now on that the Newton polygon has non-empty interior.
Suppose it is as follows, where the labels on segments are the
negative of their slopes:
$$\let\o\circ
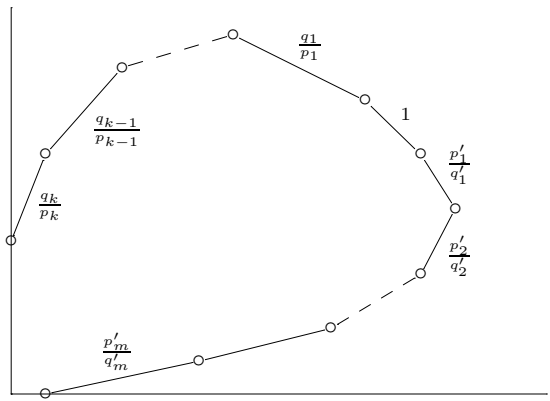
\xymatrix@=8pt@M=0pt@W=0pt@H=0pt{
\lineto[ddddddddddddd]
 &&&&&&&&&&&&&&&&&&&&&\\
&&&&&&&&\o\lineto[rrrrrdd]^{\frac{q_1}{p_1}}&&&&&&&&&&&&&\\
&&&&\o\dashto[urrrr]&&&&&&&&&&&&&&&&&\\
&&&&&&&&&&&&&\o\lineto[rrdd]^{1}&&&&&&&&\\
&&&&&&&&&&&&&&&&&&&&&\\
&\o\lineto[uuurrr]_{\frac{q_{k-1}}{p_{k-1}}}
&&&&&&&&&&&&&&\o\lineto[rdd]^{\frac{p'_1}{q'_1}}&&&&&&\\
&&&&&&&&&&&&&&&&&&&&&\\
&&&&&&&&&&&&&&&&\o\lineto[ddl]^{\frac{p'_2}{q'_2}}&&&&&\\
\o\lineto[uuur]_{\frac{q_k}{p_k}}&&&&&&&&&&&&&&&&&&&&&\\
&&&&&&&&&&&&&&&\o\dashto[ddlll]&&&&&&\\
&&&&&&&&&&&&&&&&&&&&&\\
&&&&&&&&&&&&\o\lineto[dlllll]&&&&&&&&&\\
&&&&&&&\o\lineto[dllllll]_{\frac{p'_{m}}{q'_{m}}}&&&&&&&&&&&&&&\\
\lineto[rrrrrrrrrrrrrrrrrrrr]
 &\o&&&&&&&&&&&&&&&&&&&&\\
}
$$

\smallskip Then part of (a possibly unreduced version of) the rooted
splice diagram is as follows. This diagram can be reduced if $p_k$ or
$p'_m$ is $1$.
$$\let\o\circ
\def\dec{\dotto[dl]!<6pt,-6pt>\dotto[dr]!<-6pt,-6pt>}
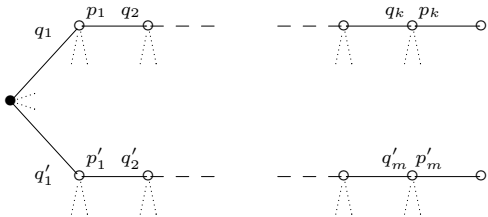
\xymatrix@R=6pt@C=8pt@M=0pt@W=0pt@H=0pt{
  &&&\Circ\lineto[rrr]^(.25){p_1}^(.75){q_2}\dec
  &&&\Circ\dashto[rrr]\dec
  &&&&&&\dashto[rrr]&&&
  \Circ\lineto[rrr]^(.75){q_k}\dec&&&
  \Circ\lineto[rrr]^(.25){p_k}\dec&&&\Circ\\ &&&&&&&&&&&&&&&&&&&\\ 
&&&\\&&&&&&&&&&&&&\\
  \Dot\lineto[uuuurrr]^(.75){q_1}\lineto[ddddrrr]_(.75){q'_1}
  \dotto[r]!<3pt,4pt>\dotto[r]!<3pt,-4pt>&&\\ &&&&&&&&&&&&&&&\\ 
&&&\\ &&&\\
  &&&\Circ\lineto[rrr]^(.25){p'_1}^(.75){q'_2}\dec
&&&\Circ\dashto[rrr]\dec&&&&&&\dashto[rrr]&&&
  \Circ\lineto[rrr]^(.75){q'_m}\dec&&&\Circ\lineto[rrr]^(.25){p'_m}\dec
&&&\Circ
\\&&&&&&&&&&&&&&&&&&& \\ \\ }$$

The number of additional edges (indicated by dotted lines) at each
vertex is computed as follows. Consider a vertex with adjacent weights
$p,q$, say. It corresponds to an edge of the Newton polygon of slope
$-q/p$ or $-p/q$.  The monomials of $f(x,y)$ for this edge sum to an
expression of the form
$$cx^ay^b\prod_{k=1}^s(x^p+c_ky^q)^{r_k}\quad\text{or}\quad
cx^ay^b\prod_{k=1}^s(x^q+c_ky^p)^{r_k}.$$
Then $s$ is the number of
edges in question.

The segment of slope $-1$ may have length $0$. And it will occur at an
end of the chain of edges shown (rather than in the middle) if all the
slopes shown are greater than $-1$ or all are less than $-1$, in which
case the root vertex has valency $1$.  For example, if $p,q$ are
coprime integers with $q<p$ and $f(x,y)$ has Newton polygon
$$\let\o\circ
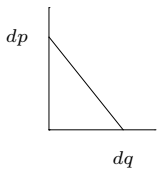
\xymatrix@=8pt@M=0pt@W=0pt@H=0pt{
&\lineto[ddddd]
\\
\scriptstyle dp& \lineto[ddddrrr]
\\ \\ \\ \\
&\lineto[rrrr]
  &&&&&\\
  &&&&\scriptstyle dq\\
}
$$
 then
$$\rsplice(f)= \def\dec{\dotto[dl]!<6pt,-6pt>\dotto[dr]!<-6pt,-6pt>}
\xymatrix@R=6pt@C=8pt@M=0pt@W=0pt@H=0pt{\Dot\lineto[rrr]^(.75){q}&&&
\Circ\dec\lineto[rrr]^(.25){p}&&&\Circ\\&&&&}~~.$$

\smallskip\noindent This is the type of polynomial studied in 
\cite{shpilrain-yu-99}.
\smallskip

\section{Algorithms}

\subsection{Linear positioning}\label{subsec:normalisation}

We first perform a linear automorphism to put points at
infinity of $f$ in ``standard position''.

Write $f(x,y)=f_0(x,y)+f_1(x,y)+\dots+f_N(x,y)$ where $f_i(x,y)$ is
the homogeneous part of degree $i$.  The points at infinity for any
fiber $f^{-1}(c)$ are the points of vanishing of the linear factors of
$f_N(x,y)$.  By a linear change of coordinates one of them can be put
at $[x{:}y]=[1{:}0]$ and, if there is more than one, another at
$[x{:}y]=[0{:}1]$.  For the Newton polygon this means that the part
lying on the line $x+y=N$ either does not include the points $(0,N)$
and $(N,0)$ or consists only of the former point.

We shall assume from now on that $f$ has been positioned in this way.

\subsection{Reduction of degree}\label{subsec:reddeg}
A reduced rooted splice diagram is \emph{minimal} if it results from
the reduced unrooted splice diagram by placing a root vertex either at
an existing vertex or on an existing edge.  We also call the
exceptional diagrams discussed in Subsection \ref{subsec:exceptional}
minimal.

\begin{proposition} Assume we have applied a linear automorphism as
  above so $[1{:}0]$ is a point at infinity for $f(x,y)$.  Then, if
  the reduced rooted splice diagram for $f(x,y)$ is not minimal, there
  is an automorphism of the form $(x,y)\mapsto
  (x-cy^q,y)$ which reduces the degree of $f$. 
\end{proposition}
\begin{proof}
  If the reduced rooted splice diagram is not minimal, the root vertex
  has just one edge emanating and that edge has far weight $1$.  In
  terms of the Newton polygon this means that for some $q>1$ there is
  a segment of slope
  $-q$ touching the $y$-axis ($q$ is the near weight at the vertex
  adjacent to the root).
  
  Then $f(x,y)$ has the form
  $$f(x,y)=by^a\prod_{k=1}^s(x+c_ky^q)^{r_k}+
  \sum_{qi+j<qS+a}c_{ij}x^iy^j,\quad\text{where } S=\sum_{k=1}^sr_k.$$
  Composing $f$ with any of the automorphisms
  $(x,y)\mapsto(x-c_ky^q,y)$ will reduce the degree of $f$.
\end{proof}

By iterating the procedure of the above proposition we can apply
automorphisms to $f(x,y)$ until its reduced rooted splice diagram is
minimal.  This procedure reduces the degree of $f$.  Thus
some minimal rooted splice diagram corresponds to an automorphic image
of $f$ of least possible degree.  However, there may be several
minimal rooted diagrams, and not all of them will give the absolute
minimal degree among automorphic images of $f$.

In any case, we assume from now on that $\rsplice(f)$ is minimal.

\subsection{Changing reduced splice diagram}
A reduced unrooted splice diagram $\splice=\splice(f^{-1}(c))$ may
have several places where one can put a root vertex to make it into a
rooted splice diagram. In \cite{neumann1} it is shown how to realise
each resulting rooted splice diagram as
$\rsplice\bigl((f\phi)^{-1}(c)\bigr)$ for a suitable polynomial
automorphism $\phi$. Since this is an essential part of our algorithm
we describe it here.

We will assume $\splice$ is non-exceptional since each exceptional
$\splice$ has a unique minimal rooted splice diagram (see Subsection
\ref{subsec:exceptional}).
  
We first describe how to find all potential root vertices for
$\splice$. Choose one possible root vertex. Consider a simple path
from this vertex to another vertex with the properties:
\begin{itemize}
\item The edge weights \emph{on} the path are all positive except
  maybe at the very end of the path;
\item If an edge weight \emph{adjacent} to the path is $>1$ then it is
  the only edge weight $>1$ at its node. A node where this happens is
  called a \emph{contributing node} (following \cite{neumann1}) and we
  say the path ``passes'' this contributing node. (If a path through a
  contributing node has the weight $>1$ on it rather than adjacent to
  it, we say it ``goes through'' the contributing node.)
\end{itemize}

Take the subgraph $\splice_0$ of $\splice$ that is the union of all
these paths.  Each possible position for a root node lies at a vertex
or on an edge of this subgraph.  Their positions are found as follows.
A \emph{directed path} in $\splice_0$ is a simple path starting at a
vertex or edge midpoint and departing each node along it by an edge
with maximal weight among the edge-weights of $\splice$ at that node.
The condition of negative edge determinants implies that a maximal
directed path will either start at a leaf, in the middle of an edge,
or at a node with all adjacent weights $1$ and it will end at a vertex
with non-positive edge weight, an arrowhead, or at a contributing
vertex.  The starting points of maximal directed paths in $\splice_0$
are the possible positions for root vertices. In fact every half-edge
of $\splice_0$ is on a unique maximal directed path, so these maximal
directed paths give a flow on $\splice_0$ and the possible root
vertices are the sources of this flow (see \cite{neumann1}).

Any simple path in $\splice_0$ will be a union of directed paths laid
end-to-end with directions alternating. The direction can only change
from forward to backward as we pass a contributing vertex and can only
change from backward to forward as we go through a potential root
vertex.  We shall call two potential root vertices ``adjacent'' if the
path between them goes through no other potential root vertex, so it
changes direction exactly once, and this happens as it passes some
contributing node. (The path may, however go through other
contributing nodes without ``passing'' them---see Item 2.\ above.)

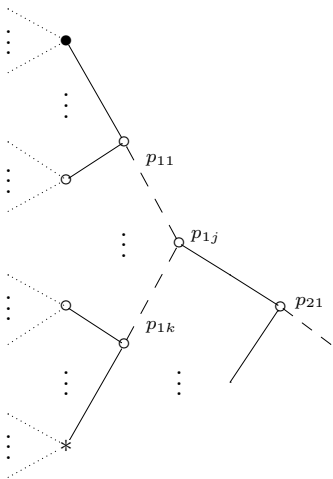
\begin{figure}[htbp]
$$
\xymatrix@R=2pt@C=18pt@M=0pt@W=0pt@H=0pt{
&\\
\Vdots&\Dot\dotto[lu]\dotto[ld]\lineto[rddd]\\ 
&\\
&\Vdots\\
&&\Circ\dashto[dddr]^(.25){p_{11}}&&&\\
\Vdots&{\Circ}\dotto[lu]\dotto[ld]\lineto[ru]\\
\\
&&\Vdots&\Circ\lineto[rd]^(.25){p_{1j}}&\\
&&&&\lineto[rd]\\
\Vdots&{\Circ}\dotto[lu]\dotto[ld]\lineto[rd]&&&&\Circ\dashto[rd]^(.25){p_{21}}
\lineto[ldd]\\
&&\Circ\dashto[uuur]_(.25){p_{1k}}&&&&\\
&\Vdots&&\Vdots&&\\
&\\
\Vdots&{*}\dotto[lu]\dotto[ld]\lineto[ruuu]\\
\\}
$$
\caption{Adjacent potential root vertices}\label{fig:roots}
\end{figure}

A typical path between adjacent potential root vertices is the path
from $\bullet$ to $*$ in Figure \ref{fig:roots} (the $p_{ij}$
are greater than $1$ and increase to the right in the graph). We are
omitting vertices that are neither contributing nodes nor potential
root vertices in this picture (they are all of valency at most $2$ in
$\splice_0$). The vertices at the left of the picture are pairwise
adjacent potential root vertices.

If $f$ does have several minimal rooted splice diagrams we need to
find automorphisms that move us from one to another.  It suffices to
move between adjacent potential root vertices.

Suppose the root of $\rsplice(f)$ is at the top left in Figure
\ref{fig:roots}. Let the rightward edge at that node correspond to the
point $[1{:}0]$ at infinity.  Then the Newton polygon will have
segments on its boundary with slopes $-p_{11}, \dots,-{p_{1j}},\dots,
-p_{21},\dots$.  There may be segments with other slopes corresponding
to vertices of $\rsplice(f)$ that are valency two in $\splice_0$ and
therefore not drawn in the figure or a segment of slope 1 coming from
additional edges at vertex $\bullet$.  Suppose we want to move the root
vertex to the position $*$. This moves it past the contributing node
$v$ with adjacent weight $p_{1j}$. The terms of $f(x,y)$ corresponding
to the Newton polygon segment of slope $-p_{1j}$ give a polynomial
$x^ay^b\prod_{i=1}^s(x+c_iy^{p_{1j}})^{r_i}$ with each factor
$(x+c_iy^{p_{1j}})^{r_i}$ corresponding to one of the relevant edges
emerging to the left from $v$.  If the edge at $v$ towards $*$ is the
one with index $i$ then the automorphism
$(x,y)\mapsto(x-c_iy^{p_{1j}},y)$ moves the root vertex as desired.
However, the resulting diagram may need further automorphisms as in
Subsection \ref{subsec:reddeg} to reduce it to a minimal diagram again
(an automorphism that gives a minimal diagram will be of the form
$(x,y)\mapsto (x-g(y),y)$ with $g(y)$ a polynomial with highest order
term $c_iy^{p_{1j}}$).  Moreover, the result may not be to move the
root vertex to $*$, but rather to some other of the potential roots to
the left of $v$ via the $i$-th edge out of vertex $v$, for instance,
the vertex above $*$ in the figure.  But we can iterate the procedure
to try to move from the new root to $*$.  Since each iteration
involves a contributing vertex closer to $*$, this iterative procedure
eventually ends with our node at $*$.

\subsection{Comparing polynomials}
Given two polynomials $f$ and $g$ in $\C[x,y]$, their (reduced
unrooted) splice diagrams are likely to distinguish them if they are
not equivalent (or even not right-left-equivalent), so the first thing
to do is compute these and reduce them as above.  Often $f$ and $g$
can be distinguished without computing the whole splice diagram, or
they may be distinguished by an irregular splice diagram or by
positions of singularities or of irregular fibers.  We will just
consider the regular splice diagram here.

We may thus assume that we have applied automorphisms to $f$ and $g$
so they have isomorphic rooted regular splice diagrams and these
diagrams are minimal.  After fixing $g$, there may be several ways of
introducing a root in the unrooted splice diagram for $f$ to make the
diagrams isomorphic, and several isomorphisms, so there may be several
cases to try.  For each case we proceed as follows.
 
First, apply a linear automorphism as in subsection
\ref{subsec:normalisation}, if it exists, so that the points at
infinity of $f$ and $g$ are the same for corresponding edges at the
root vertices of $\rsplice(f)$ and $\rsplice(g)$. This automorphism
will exist if there are three or less points at infinity, but it is
uniquely determined up to scaling by three points, so if there are
four or more points at infinity it may not exist, in which case there
is no automorphism taking $f$ to $g$ (or even $f$ to $\alpha g+\beta$
with $\alpha,\beta\in\C$) for the given isomorphism of splice
diagrams.

Once points at infinity are matched, the problem of finding an
automorphism, if it exists, taking $f$ to $g$ (or $f$ to some $\alpha
g+\beta$) for this particular isomorphism of splice diagrams reduces,
by the following theorem, to solving simple equations in a very
limited number of variables ($\alpha, \beta$, and the coefficients of
the $\phi$ described in the theorem).

\begin{theorem}
  Suppose $f=f(x,y)$ has minimal rooted splice diagram.  Assume also
  that a point at infinity occurs at $[1{:}0]$, and a second, if it
  exists, at $[0{:}1]$, and that $f$ is not a polynomial in $y$ alone.
  Suppose $\phi$ is an automorphism that does not move the root vertex
  of $f$ and does not move any point at infinity\footnote{From an
    algebraic point of view, to say an automorphism moves neither the
    root vertex nor any point at infinity means that each branch at
    infinity of a generic fiber $f^{-1}(c)$ should be at the same
    point at infinity before and after applying the automorphism.}.
  Then
  \begin{itemize}
  \item If $f$ has one point at infinity then $\phi$ is in the
    Jonqui\`ere subgroup $B$, that is, it is of the form
    $\phi(x,y)=(bx+h(y),dy+t)$. The degree of $h$ is bounded by the
    absolute value of slope of the top segment of the Newton polygon
    boundary (in particular, $\operatorname{deg}(h)\le
    \operatorname{deg}(f)$).
\item If $f$ has two points at infinity then $\phi(x,y)=(ax+s,dy+t)$.
\item If $f$ has three or more points at infinity then
  $\phi(x,y)=(ax+s,ay+t)$.
  \end{itemize}
\end{theorem}
\begin{proof}
  The specific form of $\phi$ described in the theorem is easily seen
  if $\phi$ is affine, so assume it is not.  Since the automorphism
  $\phi$ does not change the splice diagram, it does not change the
  degree of $f$ (see Item \ref{it:degree} in Section
  \ref{sec:invariants}). By Theorem \ref{th:pw}, a normal form
  representation for the automorphism $\phi$ involves exactly one
  factor from $B$.  We can thus write it as
  $\phi_1\circ\phi_2\circ\phi_3$, with $\phi_2\in B-A$ and
  $\phi_1,\phi_3\in A$, but with no restriction on non-triviality of
  $\phi_1,\phi_3$. We can modify $\phi_2$ and $\phi_3$ to make
  $\phi_2$ triangular, say $(x,y)\mapsto(x+h(y),y)$ with $\deg(h)>1$.
  Denote $f_i=f\phi_1\dots\phi_i$. By Theorem \ref{th:pw},
  $f,f_1,f_2,f_3$ all have the same degree $N$ say.  Note that $f_1$
  has the same number of points at infinity as $f$. If $f_1(x,y)$ had
  a point at infinity other than $[1{:}0]$ then it would have a
  monomial $x^ay^{N-a}$ with $a>1$. Taking the largest such $a$, it
  would follow that $f_2=f_1\phi_2$ has a monomial $y^{a\deg(h)+N-a}$
  and thus has degree larger than $N$.  Thus $f_1$ has just one point
  at infinity at $[1{:}0]$, so $\phi_1$ fixes this point and is
  therefore in $A\cap B$. The same argument applied to
  $f=f_3\psi_3\psi_2\psi_1$ with $\psi_i=\phi_i^{-1}$ shows $\psi_3\in
  A\cap B$ so $\phi_3\in A\cap B$.  Thus $\phi\in B$, as claimed.
  
  The statement about the degree of $h(y)$ follows by noting that, if
  the slope of the top segment in the Newton polygon is $-p/q$ and if
  $\deg(h)> p/q$ then a monomial $x^{N-rq}y^{rp}$ corresponding to the
  lowest point on this top segment of the Newton polygon for $f$ leads
  to a monomial $y^{(N-rq)\deg(h)+rp}$ for $f\phi$. Since
  $(N-rq)\deg(h)+rp>(N-rq)p/q+rp=Np/q>N$, this completes the proof.
\end{proof}
\section{Examples and discussion}

\subsection{Uniquely one point at infinity}
If $p,q$ are coprime integers with $1<q<p$ and $f(x,y)$ has Newton
polygon
$$\let\o\circ
\xymatrix@=8pt@M=0pt@W=0pt@H=0pt{
&\lineto[ddddd]
\\
\scriptstyle dp& \lineto[ddddrrr]
\\ \\ \\ \\
&\lineto[rrrr]
  &&&&&\\
  &&&&\scriptstyle dq\\
}
$$
then
$$\rsplice(f)= \def\dec{\dotto[dl]!<6pt,-6pt>\dotto[dr]!<-6pt,-6pt>}
\xymatrix@R=6pt@C=8pt@M=0pt@W=0pt@H=0pt{\Dot\lineto[rrr]^(.75){q}&&&
  \Circ\dec\lineto[rrr]^(.25){p}&&&\Circ\\&&&&}~~$$
is the unique
minimal splice diagram, so $f$ cannot be equivalent (or even
right-left-equivalent) to another polynomial of the same type with
different $p$, $q$, or $d$. This is the case proved in
\cite{shpilrain-yu-99}.  

More generally, if the top segment of the Newton polygon boundary has
slope $-p/q$ with $p,q$ as above then there is one point at infinity
and the minimal splice diagram is unique. Thus the invariants
extractable from this diagram are invariants of the equivalence class
of $f$. In particular, $p$, $q$, $\deg(f)$, and the lines of negative
slope on the Newton polygon boundary are invariants.

\subsection{Recognising a coordinate}
If $g(x,y)=x$ then
$\rsplice(g)=\xymatrix@1@M=0pt@W=0pt{\bullet\ar[r]&}$~,
 so if
$f$ is equivalent to $g$ and $\operatorname{deg}(f)>1$ then, after
applying a linear automorphism as in Subsection
\ref{subsec:normalisation}, there must be an automorphism that reduces
degree derivable from the Newton polygon as in Subsection
\ref{subsec:reddeg}.

The Abhyankar-Moh-Suzuki theorem (\cite{abhyankar-moh}, \cite{suzuki})
says that if $f\in\C[x,y]$ has a non-singular contractible fibre then
$f$ is a coordinate (meaning equivalent to $x$).  There are now many
proofs of this, and they all proceed by (implicitly or explicitly)
constructing the linear and triangular automorphism which together
reduce the degree of $f$. A rather different algorithm to recognise a
coordinate is given in \cite{shpilrain-yu-97}.

\subsection{Isomorphic non-equivalent curves with one place at
  infinity} In \cite{shpilrain-yu-99b} Shpilrain and Yu give the
following interesting family of examples.  Let
$$f(x,y)=x-h(x^q,y)\quad\text{and}\quad g(x,y)=x-h^q(x,y).$$
Then $f$
and $g$ are isomorphic. Recall this means the algebraic sets
$f^{-1}(0)$ and $g^{-1}(0)$ are isomorphic as affine schemes, that is,
there is a ring isomorphism
$$\C[x,y]/\bigl(x-h(x^q,y)\bigr)~\to~\C[x,y]/\bigl(x-h^q(x,y)\bigr).$$
It is given by 
$$
\begin{array}{l}
x\mapsto h(x,y)\\ y\mapsto y   
\end{array}
\quad\text{with inverse}\quad
\begin{array}{l}
  x\mapsto
x^q\\ y\mapsto y
\end{array}
$$  \WDN It would be nice to see how splice diagram
behaves in general for these.

In the earlier paper \cite{shpilrain-yu-99} they applied this in the
special case
$$f(x,y)=x-x^{q_1q_2\dots q_n}-y^p,\quad gcd(q_1q_2\dots q_n,p)=1$$
and $q=q_1q_2\dots q_i$ for $i=1,2,\dots,n-1$, to give examples of
arbitrarily many affine plane curves with one place at infinity that
are isomorphic but differently embedded in $\C^2$ (by Abhyankar and
Singh \cite{abhyankar-singh} the number of embeddings of a one-place
curve is always finite).

The rooted splice diagrams for this example are, with $Q=q_1\dots q_n$,
$q=q_1\dots q_i$:
\begin{align*}\rsplice\left(x-(x^q+y^p)^{Q/q}\right)&=
\xymatrix@R=24pt@C=24pt@M=0pt@W=0pt@H=0pt{
  {*}\lineto[r]^(.75){q}&\Circ\lineto[r]^(.75){p}\lineto[d]^(.3){p}&
  \Circ\ar[r]\lineto[d]^(.3){Q/q}&\\
  &{*}&\Circ}\\
&=\xymatrix@R=24pt@C=24pt@M=0pt@W=0pt@H=0pt{
  {*}\lineto[r]^(.75){Q}&
  \Circ\ar[r]\lineto[d]^(.3){p}&\\
  &{*}}\qquad\text{if~~~}i=n,
\end{align*}
with root vertex at whichever of the two vertices marked $*$ has the
smaller adjacent weight.

\subsection{Families of isomorphic non-equivalent curves}

In contrast to the situation for one place at infinity, a curve with
several places at infinity may have a family of inequivalent
embeddings in $\C^2$.  An example is given by $\C-\{0,1\}$ embedded in
$\C^2$ as the fiber $f^{-1}(c)$ of $f(x,y)=(xy+1)(x(xy+1)+1)$ for any
$c\ne0,1$.

Indeed, one can check that an explicit isomorphism
$f^{-1}(c)\to\C-\{0,1\}$ is $(x,y)\mapsto z:=(xy+1)/c$.  The
generators $z,z^{-1}$, and $(1-z)^{-1}$ of the ring of functions on
$\C-\{0,1\}$ are given in terms of the functions $x,y$ on $f^{-1}(c)$,
and conversely, by
\begin{align*}
z&=\frac{xy+1}{c}
\\z^{-1}&=x(xy+1)+1
\\(1-z)^{-1}&=\frac{c^2+cxy+y+xy(xy+1)}{c(c-1)}~,
\end{align*}\begin{align*}
x&=\frac1c(1-z)z^{-2}\\
y&=cz^2(cz-1)({1-z})^{-1}.
\end{align*}
Composing and simplifying, one finds an isomorphism
$$\C[x,y]/(f(x,y)-c)\cong \C[X,Y]/(f(X,Y)-d)$$ given by
$$X=\frac {cx}d ,\qquad
Y=\frac{cd(d-1)y+d(d-c)(xy+1)^2}{c^2(c-1)}.$$

Different fibers $f^{-1}(c)$ and $f^{-1}(d)$ are not even
right-left equivalent, since the map $\psi$ of a right-left
equivalence
\begin{equation*}
  \xymatrix{
\C^2\ar[r]^{\phi}\ar[d]^f&\C^2\ar[d]^f \\
\C\ar[r]^\psi&\C
}
\end{equation*}
would have to fix the two irregular values $0$ and $1$ of $f$ and
would thus be the identity, so it cannot map $c$ to $d$ if $c\ne d$.

This polynomial $f(x,y)$ has regular splice diagram
$$\rsplice(f)=\xymatrix@R=24pt@C=24pt@M=0pt@W=0pt@H=0pt{
  \Circ\lineto[r]^(.75){2}&\Circ\ar[d]\lineto[r]^(.25){-3}&\Dot\lineto[r]
  ^(.75){-1}
  &\Circ\ar[d]\ar[r]^(.25){2}&\\
  &&&\\}$$
and a single irregular diagram at the smooth reducible
fiber $f^{-1}(0)$:
$$\rsplice(f^{-1}(0))=\xymatrix@R=24pt@C=24pt@M=0pt@W=0pt@H=0pt{
  &\ar[l]\Circ\ar[d]\lineto[r]^(.25){-2}&\Dot\lineto[r] ^(.75){-1}
  &\Circ\ar[d]\ar[r]^(.25){2}&\\
  &&&\\}$$
The fiber $f^{-1}(1)$ is also reducible, with a normal
crossing singularity at $(0,-1)$, but it is regular at infinity.  The
contributions from this singularity and the singularity at infinity of
the fiber over $0$ fulfil Suzuki's formula (end of Section
\ref{sec:invariants}), so we see again that there are no other
non-generic fibers.

This is just one example of a polynomial with all generic fibers
isomorphic. Such polynomials are called \emph{isotrivial} and have
been classified (see Kaliman \cite{kaliman}).

\subsection{Other fields}
Our base field has been $\C$ for convenience and to emphasize the
geometric underpinnings for the approach.  However, algebraic closure
is certainly not essential since one can always perform computation in
extension fields as necessary. The assumption of characteristic zero
is also inessential for the algorithm: the proof of Theorem
\ref{th:pw} in \cite{wightwick} goes through with no change and the
algorithm, as described there without splice diagrams, also works.
However, the theory of splice diagrams would need some modification in
finite characteristic.

\end{document}